\documentclass{amsart}
\usepackage{amsmath,amsthm, amscd}
\usepackage{marvosym}
\usepackage{graphicx}
\usepackage{color}
\usepackage{epsfig}
\usepackage{fullpage}
\theoremstyle{plain}
\newtheorem{thm}{Theorem}[section]
\newtheorem{cor}[thm]{Corollary}
\newtheorem{defi}[thm]{Definition}

\newtheorem{prop}[thm]{Proposition}
\newtheorem{lemma}[thm]{Lemma}

\newtheorem*{claim*}{Claim}

\theoremstyle{definition}
\newtheorem{defn}[thm]{Definition}

\newcommand{\comment}[1]{}

\usepackage{citesort}

\begin{document}
\title{Tangle solutions for composite knots: application to Hin recombination}
\author[D.~Buck]{Dorothy~Buck}
\address{Department of Mathematics, Imperial College London}
\email{d.buck@imperial.ac.uk}
\author[D.~Mauricio]{Mauro Mauricio}
\address{Department of Mathematics, Imperial College London}
\email{mauro.mauricio@imperial.ac.uk}

\begin{abstract}
We extend the tangle model, originally developed by Ernst and Sumners \cite{Sum}, to include composite knots. We show that, for any prime tangle, there are no rational tangle attachments of distance greater than one that first yield a 4-plat and then a connected sum of 4-plats. This is done by building on results on exceptional Dehn fillings at maximal distance. We then apply our results to the action of the Hin recombinase on mutated sites. In particular, after solving the tangle equations for processive recombination, we use our work to give a complete set of solutions to the tangle equations modelling distributive recombination.
\end{abstract}

\maketitle
\today

\section{Introduction}
The tangle model was introduced by Ernst and Sumners \cite{Sum} as a tool for understanding the action of certain proteins, \textbf{site-specific recombinases}, on DNA. If the initial DNA molecule is circular, this reaction can form DNA knots and links. In this setting, the circular DNA molecules are decomposed into a collection of tangles, and the protein action is modelled as replacing one particular tangle $P$ with another, $R$. Information on the knot or link type of the resulting DNA can then be used in conjunction with the tangle model to help determine both structural and mechanistic information about site-specific recombination. This approach has been carried out successfully in many situations (\cite{Buck,Darcy2,Buck2,Ernst3,Cris,Vaz,Vaz2,Sum,Ernst2}).

\par These treatments have been restricted to the case where the DNA products are exclusively prime knots and links. However, experiments show that, in certain situations, it is possible to obtain products which are nontrivial connected sums \cite{Heichman}. The previous analysis relies on the Cyclic Surgery Theorem \cite{CGLS}. This approach breaks down for composite knots, essentially because their double branched cover do not have cyclic group.  In this paper, we extend the tangle model to be able to deal with these non-prime cases. In particular, we prove that if the distance between the parental tangle $P$ and the recombinant tangle $R$ is greater or equal than 2, then only the obvious locally knotted solutions may occur. We also give some restrictions on the solutions for the much more difficult distance one case. In particular, we completely describe the situation for Montesinos tangles. We illustrate this by applying the results to the Hin recombinase, a protein present in the bacteria \textit{Salmonella enterica}.

\vskip 0.1in
This work is organised as follows: in Section 2, we give an introduction to site-specific recombination, followed by a description of the Hin recombinase system, which acts as a motivation for the work to follow. In Section 3, we give a brief review of tangles, which are the natural objects of study in topological models of site-specific recombination. Section 4 contains the bulk of our theoretical work: we formulate a tangle model for the occurrence of composite products, and we relate it to the problem of exceptional Dehn fillings on simple manifolds at maximal distance. Here, we show that there are no rational tangle attachments to a prime tangle yielding a 4-plat and then a connected sum of 4-plats. In Section 5, we solve the tangle equations arising from the processive Hin recombination on mutated sites. In Section 6, we combine this with our results from Section 4 to obtain the complete set of tangle solutions to the equations arising from distributive Hin recombination (Theorem \ref{ti}). We conclude in Section 7 by discussing future directions for research.

\section{Biological Motivation}

\subsection{Site-specific Recombination}
Site-specific recombination is a fairly broad term that is used to describe a localised protein reaction in which two double stranded DNA segments (possibly belonging to the same DNA molecule) are brought together, cut, rearranged and resealed in some particular manner. This results in a variety of local operations on the genome of organisms, such as insertion of segments, a key step in viral infections, or inversion of a particular sequence, which can be a way to regulate gene expression (see \cite{Grindley} for more details). 

\par This process can be described in more detail as follows: in a first step, recombinases detect each individual specific DNA sequence and bind to it. Typically, four protein molecules bind to two sites (two molecules for each segment). The two segments are then brought together, either by protein action or random thermal motion, and juxtaposed. The juxtaposed sites are known as the \textbf{crossover sites}. The trapped DNA together with the protein is known as the \textbf{synaptosome}. In a second stage of the reaction, DNA strands are cleaved, exchanged and resealed within the synaptosome complex. A total of four such reactions occur, one for each strand, effectively disconnecting each crossover site. A final stage of the reaction involves the breaking down of the synaptosome, and the disassociation of the protein from the DNA. Before this occurs, the protein can immediately cleave again and perform one more round of recombination. This phenomenon is known as \textbf{processive recombination}. Naturally, proteins may act on DNA which has already undergone one or more full recombination events, and this is known as \textbf{distributive recombination}.

\subsection{The Hin recombinase} Our work will focus chiefly on the particular case of Hin-mediated recombination, which we describe in more detail. Hin is a site-specific recombinase from the serine family. It is found \textit{in vivo} in members of the \textit{Salmonella} family of bacteria, for example \textit{Salmonella enterica}. Its primary function is to regulate the expression of the \textit{fliC} and \textit{fliB} genes. These genes are involved in the synthesis of the \textit{flagellum} of the bacteria. The \textit{flagellum} is a tail-like structure present in some bacteria, and it is their primary source of locomotion through media. A primary constituent of the \textit{flagellum} is the protein known as flagellin. This helps the \textit{flagellum} achieve a helicoidal structure, which favours locomotion. It is precisely the synthesis of this protein that fliC and fliB regulate.

\par Depending on the state of activity of these genes, slightly different versions of flagellin will be assembled. This phenomenon is known as \textbf{flagellar phase variation}. Recombination by Hin inverts a particular DNA sequence, thereby acting as an ``on-off'' switch for gene expression. It acts on a so-called $H$ segment, containing the regions that code for \textit{fliB} and \textit{fliC}. It is known that, when $H$ is in standard position, expression of \textit{fliC} is inhibited, resulting in one type of flagellar phase, whereas when $H$ is in inverted position, \textit{fliB} is inhibited and \textit{fliC} is expressed, resulting in the alternative flagellar variety. This is described in more detail by Kutsukake \textit{et al} in \cite{Katsu}.

\vskip 0.1in
 $H$ is a 996 base pair (bp) sequence, with two specific sites, each of length 26bp, at its extremes. These sites are designated \textit{hixL} and \textit{hixR}, and it is known that these are the specific sites that Hin binds to. Figure 1 illustrates the known bp sequence for \textit{hixL}. The two letter subsequence $AA$ (or $TT$, depending on the strand) in bold is the region where DNA cleavage takes place. Cleavage occurs along the region specified by the gray line. The \textit{hixR} site has a slightly different sequence, but it also contains the cleavage sites highlighted.

\begin{figure}[h!]
\centering
\includegraphics[scale=0.6]{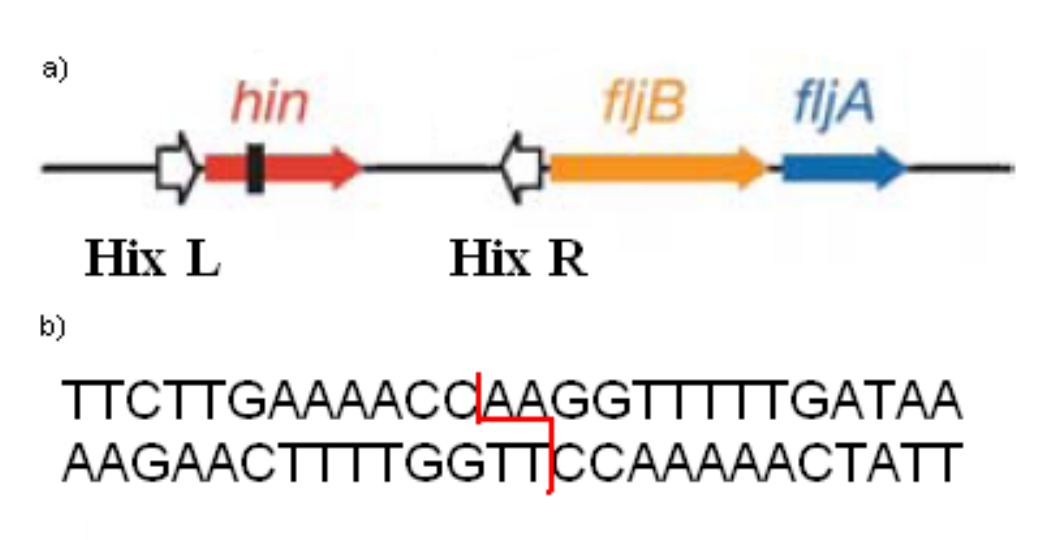}
\caption{a) Schematic of the \textit{H} sequence, together with the \textit{hixL} and \textit{hixR} sequence. The \textit{hixL} sequence. b) The \textit{hixL} sequence spelt in detail, together with the cleavage sites in gray. Adapted from \cite{Grindley}.}
\label{hixl}
\end{figure}

\par The Hin mediated inversion reaction is catalysed by enhancer proteins HU and Fis. In the $H$ segment, there is also a 65bp long \textbf{enhancer sequence}. During the inversion reaction, two pairs of molecules of Fis bind to the enhancer sequence to help form the synaptic complex: see Figure 3 (the enhancer sequence is unaffected during subsequent stages of recombination). The synaptosome consists then of four subunits of Hin, bound in pairs to each of the recombination sites. Following the terminology of \cite{reid}, we will henceforth name it the \textbf{invertasome}.

\begin{figure}[h!]
\centering
\includegraphics[scale=.8]{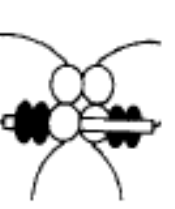}
\caption{Proposed model for the configuration of the invertasome. The white balls are the four subunits of Hin, whereas the solid black balls are units of Fis bound to the enhancer sequence (depicted by the white rectangle).}
\label{inve}
\end{figure}

\subsection{Hin-mediated recombination yields composite DNA knots}
 
\begin{figure}[h!]
\centering
\includegraphics[scale=0.6]{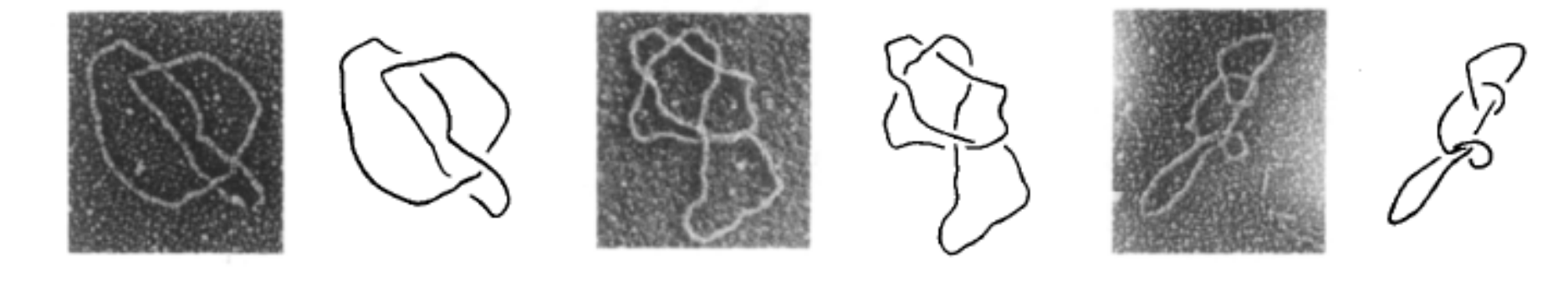}
\caption{Rec-A enhanced electron microscopy images of the products of the pRJ862 reaction. One sees a trefoil knot, the $5_2$ twist knot, and a composite knot $3_1\#3_1$. Taken from \cite{Heichman}}
\label{lhtrefoil}
\end{figure}

\vskip 0.1in In \cite{Heichman}, Heichman \textit{et al} performed a series of experiments with Hin \textit{in vitro}, on both ``standard'' and mutated circular DNA molecules. To create the mutants, Heichman and colleagues introduced a small mutation at the cleavage site of \textit{hixL}, replacing $AA$ by $AT$. These molecules were then recombined by Hin, and the resulting knot types were analysed both by agarose gel techniques and by electron microscopy\footnote{It was also observed that the standard type recombination yields a sequence of products of the form $0_1\rightarrow3_1\rightarrow4_1\rightarrow5_2$. The tangle analysis of this system is essentially carried out in Vazquez and Sumners \cite{Vaz}, and we refer the reader to it for details.}. The electron microscopy results for recombination on mutated sites are displayed in Figure \ref{lhtrefoil}. Therefore, we have a sequence of reaction products as $0_1\rightarrow3_1\rightarrow5_2$. There is also the occurrence of a non-prime knot $3_1\#3_1$. The agarose gels also reveal the formation of a 7 crossing knot, the knot type of which was not imaged via electron microscopy.  

\begin{figure}[h!]
\centering
\includegraphics[scale=0.6]{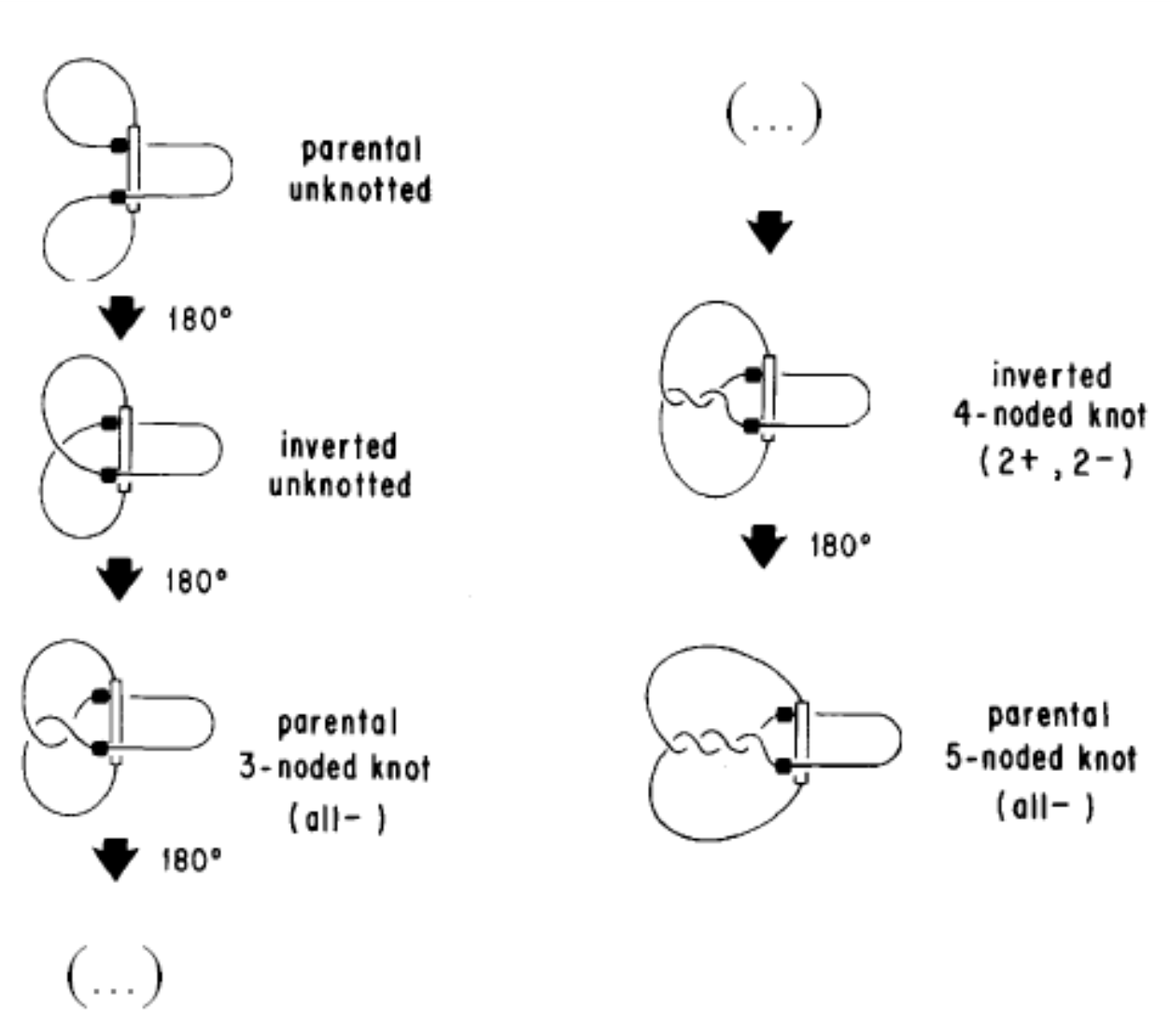}
\caption{The model proposed in \cite{Heichman} predicts a family of products in accordance with the observed data.}
\label{invertassome}
\end{figure}

\begin{figure}[h!]
\centering
\includegraphics[scale=0.6]{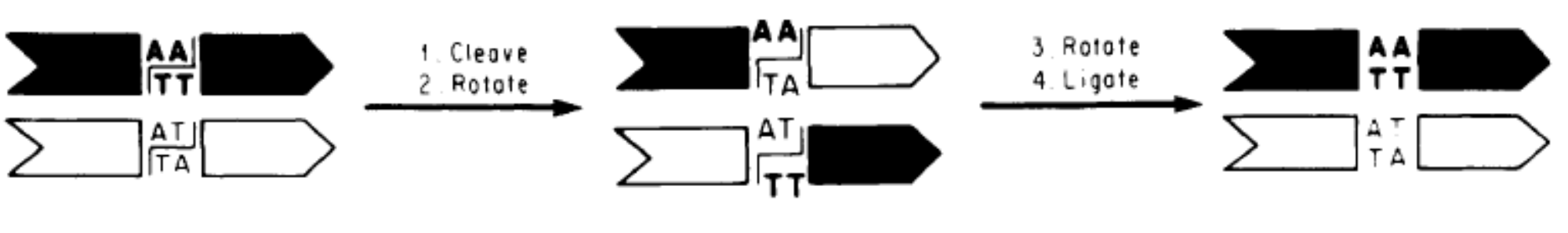}
\caption{Mismatch in the DNA sequence after one attempted round of recombination. Taken from \cite{Heichman}}
\label{rotate}
\end{figure}

The spectrum of knots observed led Heichman \textit{et al} to propose the following strand exchange mechanism for the action of Hin: cleavage occurs at the sites, in the standard manner for a serine recombinase, and then a 180$^\circ$ clockwise turn is performed at each recombination round, as illustrated in Figure \ref{invertassome}: this is in accordance with the products that were observed. With this scheme, it is also easy to explain the spectrum of products observed for the action on mutant sites, namely the absence of figure-eight knots and 6 crossing twist knots: upon performing the first round of rotation, Hin tries to reseal the strands, but it finds that there is a mismatch in the sequence (\textit{i.e.} the strands are not complementary). This causes the system to undergo an additional round of recombination in order to achieve a configuration that is consistent with base pair complementarity: see Figure \ref{rotate}.

 As for the step that produces $3_1\#3_1$ and other non-prime knots, this cannot come from processive recombination, since this only produces twist knots, which are prime. The most likely explanation is that this product arises when Hin performs distributive recombination. It is conceivable that Hin could perform one round of recombination, resulting in $3_1$, before releasing DNA. It could then bind once more to form a new invertasome, keeping the knotted arc in the outside unbound DNA. After one round of recombination, this would produce $3_1\#3_1$.

\vskip 0.1in We are thus motivated to apply the tangle model of Ernst and Sumners to give mathematical confirmation of these predictions. Our main goals are the following:

\begin{enumerate}
\item To describe the formation of the knots $3_1$, $5_2$ and the 7-crossing knot, as far as the topological mechanism of strand exchange is concerned.

\item To describe, in terms of an appropriate tangle model, the formation of composite knots $3_1\#3_1$ as reaction products.
\end{enumerate}

The first item will be in the same spirit of the aforementioned classical tangle model. The second step, however, will require a slightly new approach to the tangle model, not so much in formulation, but in the techniques required to solve it. We postpone the details of this until Section 4, after we have defined the mathematical objects under scrutiny, the tangles.

\section{Background on tangles}
Tangles were first introduced by Conway \cite{Conway}, as part of a scheme to develop a more efficient notation for enumerating knots. Since then, quite a lot of work has been done on them, and they are a rich source of examples in Dehn surgery, for they provide a fairly friendly way to visualise the properties of some 3-manifolds (for a comprehensive discussion of this, see Gordon \cite{G2}). We begin by making the definition:

\begin{defi} A \textbf{tangle} is a topological triple $(B,t,\psi)$, where $B$ is (PL or smoothly) homeomorphic to a 3-ball, $t$ is a pair of arcs properly embedded in $B$, and $\psi$ is an orientation-preserving embedding $\psi:(\partial B,\partial t)\rightarrow (S^2,P)$, where $(S^2,P)$ is the unit 2-sphere in Euclidean space, centred at the origin, together with a choice of 4 distinguished points $P=\{NW, SW, NE, SE\}$. We call $\psi$ the \textbf{boundary parametrisation} of $(B,t)$.
\end{defi}

We will commonly refer to $t$ as the \textbf{arcs} of the tangle. In situations where the boundary parametrisation is irrelevant, we will simply use the notation $(B,t)$.
We will consider two notions of equivalence. The first one, which ignores $\psi$, will be given the name of tangle homeomorphism, whereas the second (a stronger notion of equivalence) will receive the name of tangle isomorphism.

\begin{defi} Two tangles $(B_1,t_1)$ and $(B_2,t_2)$ are said to be \textbf{homeomorphic} if there exists a homeomorphism $h:(B_1,t_1)\rightarrow (B_2,t_2)$ such that $h(t_1)=t_2$.
\end{defi} 

\begin{defi} Two tangles $(B_1,t_1,\psi_1)$ and $(B_2,t_2,\psi_2)$ are said to be \textbf{isomorphic} if there exists a homeomorphism $h:(B_1,t_1)\rightarrow (B_2,t_2)$ as above with the additional requirement that the following diagram commutes ($Id$ is understood to be the identity map).
\end{defi}

$$\begin{CD}
B_1 @>h>> B_2\\
@VV\psi_1V @VV\psi_2V\\
S^2 @>Id>> S^2
\end{CD}
$$
\vskip 0.1in
Tangles are partitioned into 3 equivalence classes, up to homeomorphism:

\begin{defi} Up to tangle homeomorphism, we have the following:

\begin{enumerate}

\item  A tangle $(B,t)$ is said to be \textbf{locally knotted} if $\exists$ an embedded 2-sphere $S$ in $B$ intersecting $t$ in two points and such that $S$ bounds a (ball, knotted arc) pair. 
\item A tangle $(B,t)$ is said to be \textbf{rational} if it is homeomorphic to the tangle $(D^2\times I, \{x,y\}\times I)$. This is equivalent to the existence of a properly embedded disk $D'$ in $B$ such that $D'$ separates the arcs, and such that no 2-sphere intersects the arcs in the manner described in the previous paragraph.
\item A tangle $(B,t)$ is said to be \textbf{prime} if it is neither locally knotted nor rational.
\end{enumerate}
\end{defi}

Tangles are best studied together with their two-fold covers, branched over $t$. Notice that the contractibility of $B$ ensures that these are uniquely defined once the branching set is specified, so we can speak of \textit{the} two-fold cover, branched over $t$, without ambiguity. When the branching set downstairs is implicit in the discussion, we will abuse terminology and will refer to the cover simply as the \textbf{double branched cover} of a tangle. The following theorem of Lickorish \cite{Lick} is fundamental:

\begin{thm}\cite{Lick} Let $(B,t)$ be a tangle, and $(\tilde{B},\tilde{t})$ its double branched cover. Then, the following holds:

\begin{enumerate}
\item $(B,t)$ is locally knotted $\Leftrightarrow (\tilde{B},\tilde{t})$ is a non-prime manifold with torus boundary.
\item $(B,t)$ is rational $\Leftrightarrow (\tilde{B},\tilde{t})$ is a solid torus.
\item $(B,t)$ is prime $\Leftrightarrow(\tilde{B},\tilde{t})$ is irreducible, $\partial$-irreducible with torus boundary. (Recall that the solid torus is the only irreducible, $\partial$-reducible manifold with torus boundary.)
\end{enumerate}
\end{thm}

Rational tangles are classified up to isomorphism by the extended rational numbers $\mathbb{Q}\cup \{\infty\}$, in the following sense: we define three types of homeomorphisms of tangles

\begin{enumerate}
\item $h:(B,t,\psi_1)\rightarrow (B',t',\psi_2)$ is defined as a positive Dehn half-twist about the meridian disk $D_{\mu}$, which is the intersection of the plane $y=0$ with the standard unit 3-ball.
\item $v:(B,t,\psi_1)\rightarrow (B',t',\psi_2)$ is defined as a positive Dehn half-twist about the latitude disk $D_{\lambda}$, which is the intersection of the plane $x=0$ with the standard unit 3-ball.
\item $r:(B,t,\psi_1)\rightarrow (B',t',\psi_2)$ is defined as a reflection in the plane containing $NW$, $SE$ and the origin.
\end{enumerate}

Schematically, these are represented in Figure 6.
\begin{figure}[h!]
\centering
\includegraphics[scale=0.4]{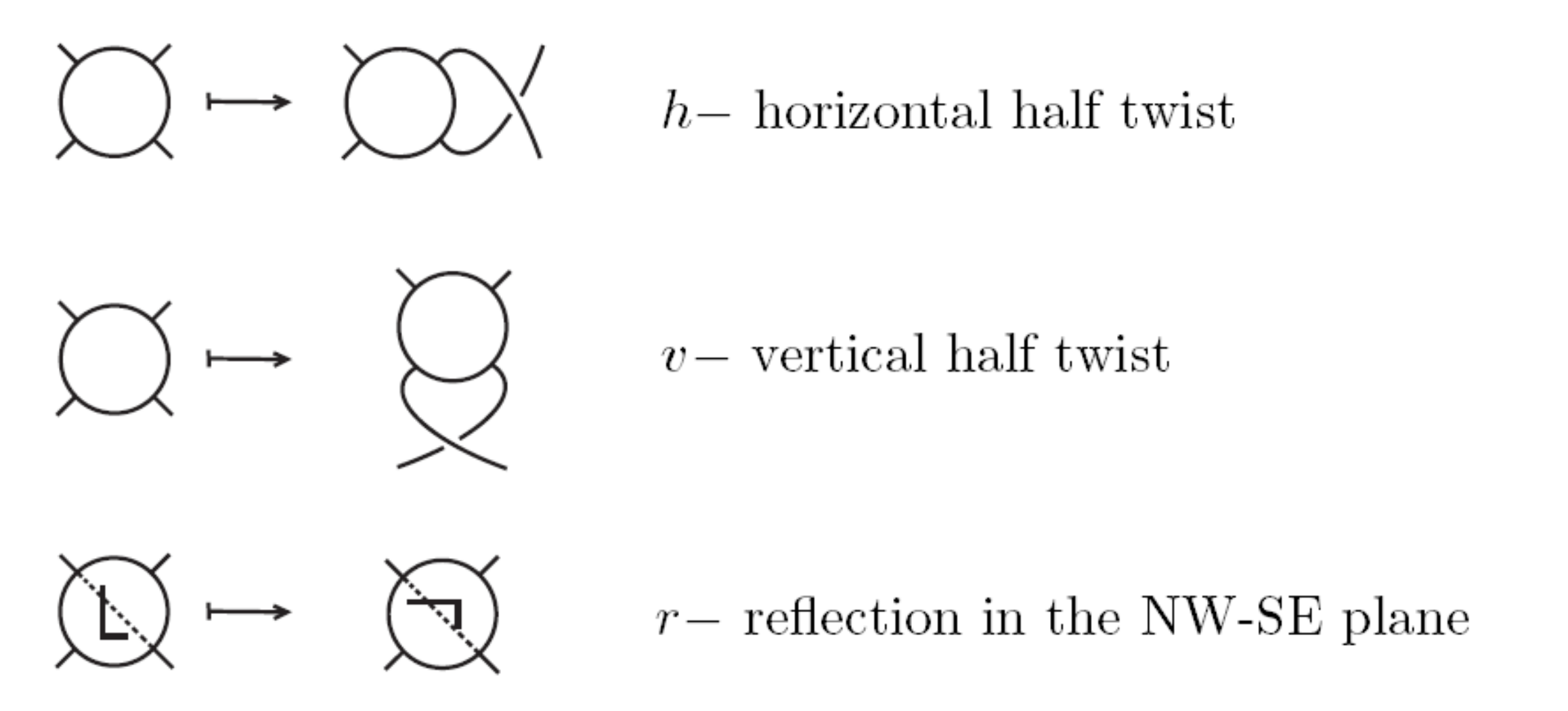}
\label{f:twists}
\caption{The twist homeomorphisms on rational tangles. Adapted from \cite{G3}.}
\end{figure}

Note that the twisting convention is that of Gordon \cite{G3}. It is opposite to that of \cite{Sum}.
\vskip 0.1in
We also define two special types of tangle, $T(0)$ and $T(\infty)$, displayed in Figure 7:

\begin{figure}[h!]
\centering
\includegraphics[scale=0.4]{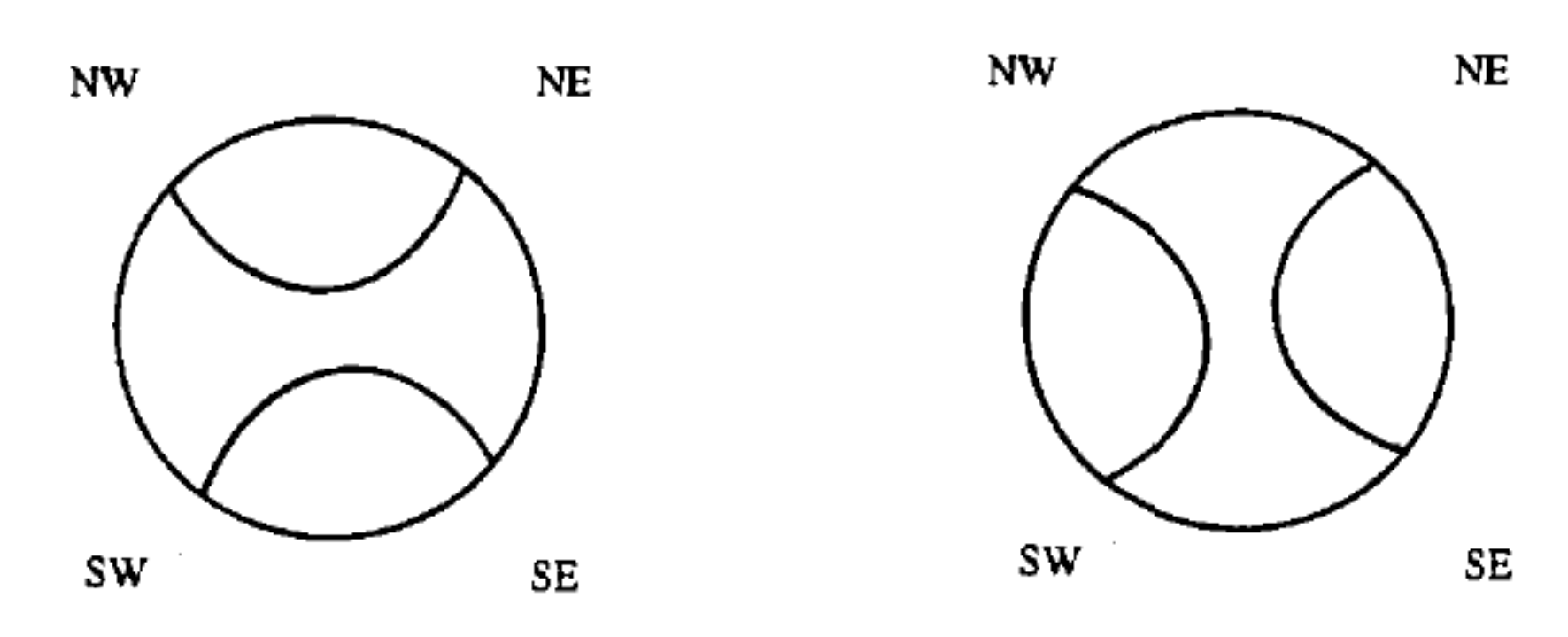}
\label{f:twists}
\caption{The tangles $T(0)$ (left) and $T(\infty)$ (right). Adapted from \cite{Sum}.}
\end{figure}

With all of this in mind, we have the following theorem of Conway \cite{Conway}:

\begin{thm} We have the following facts about rational tangles:
\begin{enumerate}
\item A tangle $(B,t,\psi)$ is rational if and only if it can be obtained from $T(\infty)$ by a finite sequence of horizontal twist homeomorphisms alternating with reflections taking the form $h^{a_1}rh^{a_2}r...h^{a_n}r$, where the $a_i\in\mathbb{Z}$. Setting $p/q=a_1+\frac{1}{a_2+\frac{1}{...+\frac{1}{a_n}}}$, with $(p,q)=1$, we write $(B,t,\psi)=T(p/q)$.
\item 
 $T(p_1/q_1)$ is isomorphic to $T(p_2/q_2)\Leftrightarrow p_1/q_1=p_2/q_2$
\end{enumerate}
\end{thm}

A rational tangle of the form $T(m/1)$ is called \textbf{integral}.

\vskip 0.1in One can define two important operations on tangles. Given tangles $A$ and $B$, we can form their \textbf{tangle sum} by connecting the $NE$ of $A$ with  the $NW$ of $B$ via a trivial arc. Similarly, one connects the $SE$ of $A$ with  the $SW$ of $B$. the result is (often) a tangle, which we denote by $A+B$. Also, given a tangle $A$, one can form its \textbf{numerator closure} by connecting the northern arcs via a trivial arc, and the southern arcs in the same fashion. The result, $N(A)$, is a topological pair $(S^3,K)$, where $K$ is some knot or link. These constructions are schematised in Figure 7.

\begin{figure}[h!]
\centering
\includegraphics[scale=0.4]{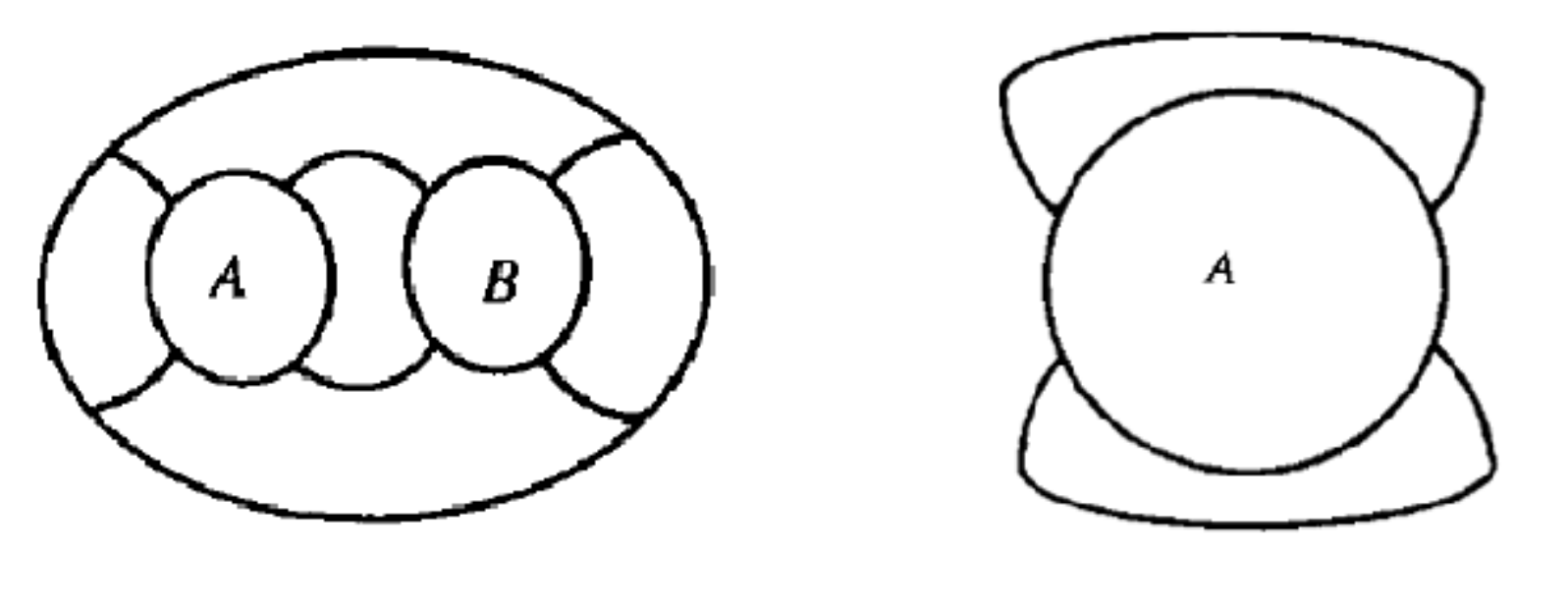}
\caption{$A+B$ (left) and $N(A)$ (right).}
\label{fig:sum}
\end{figure}

Given rational tangles $T(r_1),..,T(r_n)$, we define the \textbf{Montesinos tangle} $M(r_1,..,r_n)$ as $T(r_1)+..+T(r_n)$. It is shown in \cite{Monte} that the double branched cover $M(r_1,..,r_n)$ is a Seifert fibered space over the disk with exceptional fibers of multiplicities $1/r_1,..,1/r_n$.

We can now formulate the tangle model. We model the synaptosome as a ball, which represents the region occupied by the 4 protein molecules, with a pair of properly embedded arcs, each arc representing a double stranded segment of DNA. These form a tangle which we denote by $P$. The complement of this tangle in $S^3$, together with the DNA outside the synaptosome, form a tangle which we denote by $O$. The entire DNA segment can then be represented as the numerator closure $N(O+P)$. We model the action of the recombinase as replacing the tangle $P$ by another tangle $R$. Therefore, one round of site-specific recombination can be described by the pair of \textbf{tangle equations}

$$ N(O+P)=K_0
$$
$$ N(O+R)=K_1
$$
 
for some knots $K_0$, $K_1$. It is straightforward to extend these to $m$ rounds of processive recombination. At the $i$th recombination step ($i>1$), we think of recombination as replacing the tangle $(i-1)R$ with $iR$, where $mR$ is understood to be $\underbrace{R+...+R}_{\rm m\; times}$. Hence, we have

$$N(O+P)=K_0
$$
$$N(O+R)=K_1
$$
$$N(O+2R)=K_2
$$
$$\vdots
$$
$$N(O+mR)=K_m
$$

The purpose of the model is then to determine the tangles $P$, $R$ and $O$, given the knot type of the $K_i$. Usually $K_0=0_1$, that is, the initial DNA is taken to be unknotted.

\vskip 0.1in Following the discussion in \cite{Vaz}, we remark that experimental evidence seems to indicate that Hin behaves in a more complex way, that of a 3-string tangle (see Figure \ref{inve}). These objects are still not well understood, although there has been recent progress - e.g. see the work of Cabrera-Ibarra \cite{Ib,Ib2,Ib3} and Darcy \textit{et al} \cite{Darcy}. Nevertheless, there exists biological evidence pointing to the fact that the enhancer sequence does not play an active role in the recombination process \cite{reid}. Therefore, it is possible to perform an isotopy that removes it from the region modelling the invertasome, justifying the treatment of $P$ and $R$ as 2-string tangles.

\vskip 0.1in Solving tangle equations in full generality is a very hard problem. Fortunately, the knots usually observed in DNA recombination experiments are 4-plats (also known as 2-bridge knots), and this makes the problem more tractable. Given enough reaction steps, it is often possible to uniquely determine $O$ and $R$.

\vskip 0.1in The strategy for solving tangle equations is the following: one uses the information on knot types to argue that $P$ and $R$ are rational tangles. This means that their double branched covers are solid tori. Therefore, we convert the tangle analysis into a Dehn surgery problem. That is, we can think of attaching a solid torus $\tau^{-1}(P)$ (resp. $R$) to $\widetilde{O}$. Here, $\tau$ denotes the double branched cover map. Then, results on Dehn surgery are applied to (hopefully) conclude that $\tilde{O}$ must be a solid torus, implying that $O$ is a rational tangle. A computation, as described in \cite{Sum}, then yields the tangle solutions.

\subsection{Tangle equations for Hin action} In the light of the results of Heichman's experiment, we propose the following tangle equations for processive Hin recombination on mutated sites

\begin{equation} N(O+P)=b(1,1)
\label{1}
\end{equation}
\begin{equation} N(O+R)=b(3,1)
\label{2}
\end{equation}
\begin{equation} N(O+R+R)=b(7,3)
\label{3}
\end{equation}
\begin{equation} N(O+R+R+R)=K
\label{4}
\end{equation}

Here, we have switched to the 2-bridge notation of Schubert. Therefore, $b(1,1)$ denotes the unknot, $b(3,1)$ denotes the trefoil knot, and $b(7,3)$ denotes the $5_2$ knot. $K$ is a 7-crossing knot, and its knot will be determined as a corollary of our results. The faintness of the line in the agarose gel corresponding to 7 crossing knots (see \cite{Heichman}) supports the suggestion that $K$ is the product of 3 rounds of recombination, since one should expect a dramatic decrease in the density of products as the number of recombination rounds increases.

\vskip 0.1in 

As stated in the previous section, we will also be interested in the $b(3,1)\#b(3,1)$ case, which we postulate arises from distributive recombination. To capture the essence of this phenomenon, we propose the following tangle equations

\begin{equation}N(Q+P)=b(3,1)\label{11}
\end{equation}
\begin{equation}N(Q+R)=b(3,1)\#b(3,1)\label{22}
\end{equation}
 
Notice that $Q$ is not the same as $O$. These are two different systems, which must be treated separately. Nevertheless, the information on $P$ and $R$ is crucial for the analysis, and it is part of the motivation for considering the first system.

\vskip 0.1in

In the next section, we tackle the equations for distributive recombination in a much more general form, namely

\begin{equation}N(Q+P)=K_1\label{tng1}
\end{equation}
\begin{equation}N(Q+R)=K_2\#K_3\label{tng2}
\end{equation}

where the $K_i$ are 4-plats, and $P,R$ are rational. We postpone the solution of Equations (\ref{1})-(\ref{4}) until Section 4, where we also solve Equations (\ref{11}) and (\ref{22}) as a consequence of the work in Section 3.

\section{Solution of tangle equations yielding composite knots}
Our goal in this section is to give results on the classification of tangles $O$ which arise as solutions to the tangle equations presented below:

$$N(O+P)=K_1
$$
$$N(O+R)=K_2\#K_3
$$

where the $K_i$ are 4-plats, and $P,R$ are rational. $K_2$ and $K_3$ are both assumed to be nontrivial.
\vskip 0.1in 

\textbf{Note:} $O$ cannot be a rational tangle, since the sum of two rational tangles gives a 4-plat. So we consider $O$ to be either locally knotted or prime.

\subsection{Locally knotted solutions} We start by looking for locally knotted solutions. We begin by noting that, since (ball, knotted arc) pairs persist upon tangle summation, $K_1=K_2$ or $K_1=K_3$ for locally knotted solutions to occur. If that is the case ($K_1=K_2$, say), we can excise the common (ball, knotted arc) pairs in both equations and we are thus reduced to

$$N(O'+P)=0_1
$$
$$N(O'+R)=K_3
$$

where $0_1$ is the unknot. $O'$ is a tangle from which a (ball,knotted arc) pair has been excised. Now, we can interpret our tangle equations as Dehn surgery on a strongly invertible knot $K$ in $S^3$, the core of the solid torus given by $\tau^{-1}(P)$. Our notation for Dehn surgery on a manifold $\tilde{O}$ with a torus boundary component, along a slope $\alpha$, is $\widetilde{O}(\alpha)$. If $K_3=b(p,q)$ then taking double branched covers gives

$$\widetilde{O'}(\alpha)=S^3
$$
$$\widetilde{O'}(\beta)=L(p,q)
$$

where $\alpha,\beta$ are slopes (simple closed curves up to isotopy) in $\partial \widetilde{O'}$. Note that $O'$ is either rational or prime. We would like to conclude that $\widetilde{O'}$ is a solid torus (and hence $O'$ is rational), but this is not always true. There are numerous examples of nontrivial strongly invertible knots having lens space surgeries. Following a procedure explored by Bleiler \cite{Bl}, we can construct, given a strongly invertible projection of a knot possessing a lens space surgery, examples of prime tangles having rational tangle attachments yielding 4-plats. The reader is referred to Bleiler's work \cite{Bl} and to Gordon \cite{G2} for details on this construction.

\vskip 0.1in
 Recent results of Kronheimer \textit{et al} \cite{KMOS} show the following: if $|p|<9$, then the only non-trivial knot which possesses a $L(p,q)$ surgery is the trefoil. We note that site-specific recombination usually produces 2-bridge knots of low complexity ($|p|< 9$). This, together with a complete knowledge of surgeries on torus knots as given by Moser \cite{Mos} and the remarks in the preceding paragraph, can help one enumerate locally knotted solutions to the tangle equations (\ref{tng1}) and (\ref{tng2}) for the case $|p|<9$. Whereas it would be impractical to explicitly formulate a classification theorem, individual cases can be treated with relative ease, and the details are left to the reader.

\subsection{Prime tangle solutions} To classify prime tangle solutions we use two different strategies depending on the distance between $P$ and $R$. The \textbf{distance} between two rational tangles $P=T(p/q)$ and $R=T(r/s)$ is defined as $d(P,R)=|ps-rq|$. The distance between $P$ and $R$ is related to the distance between the Dehn filling slopes which they determine as follows:

\begin{lemma}\label{rt} Let $P$ and $R$ be rational tangles. Consider $N(O+P)$, $N(O+R)$, and the corresponding lifts $\tilde{O}(\alpha)$, $\tilde{O}(\beta)$. Then, $d(P,R)=\Delta(\alpha,\beta)$, where $\Delta(\alpha,\beta)$ is the minimal geometric intersection number between the slopes $\alpha$, $\beta\in \partial \tilde{O}$.
\end{lemma}

\begin{proof} Let us consider the boundary of a tangle $T$. The $SE-SW$ arc lifts via the branched covering map to an essential simple closed curve $k_1 \in \partial \tilde{T}$. Likewise, the $NE-SE$ arc lifts to an essential simple closed curve $k_2$. It is clear that these two curves form a basis for $H_1(\partial \tilde{T})$. Now, let us fill the boundary of $T$ with a rational tangle $p/q$. Ernst \cite{Ernst} shows that the curve that generates the kernel of the inclusion $i_*: H_1(\partial \tilde{T}(p/q)) \rightarrow H_1(\tilde{T}(p/q))$ is generated by $p k_1+ q k_2$. Thus, for $\tilde{T}_1$, $p k_1+ q k_2$ bounds a disk, and for $\tilde{T}_2$, $r k_1+ s k_2$ bounds a disk. Thus, the distance between the fillings is just the geometric intersection number of these two curves, and the result follows.
\end{proof}

From now on, we will treat the cases $d(P,R)>1$ and $d(P,R)=1$ separately.

\subsubsection{\textbf{Case $d(P,R)>1$:}}To begin, we calculate the double branched cover of a connected sum of 4-plats.
\begin{prop} Given $b(p,q)\#b(r,s)$, then it is double branch covered by $L(p,q)\#L(r,s)$. 
\end{prop}

\begin{proof} Let $L(p,q)$ be the double branched cover of $b(p,q)$. Remove a ball from $S^3$ intersecting $K_1$ in an unknotted arc. The resulting operation upstairs is to puncture $\tilde{K}_1$, since a ball double branch covers itself branched along an unknotted arc. Thus, we have that the double branched cover of a ``ball with knotted $K_1$ arc" is a punctured $\tilde{K}_1$. The same goes for $K_2$. Gluing two of these balls along their boundary gives the pair $(S^3, K_1\#K_2$), and the corresponding operation upstairs is to glue the punctured branched covers along their $S^2$ boundary. The result follows.  
\end{proof}

 Translating our problem to the double branched cover setting, we have the following Dehn Filling problem:

\begin{equation}\tilde{O}(\alpha)=L(p,q)\label{111}
\end{equation}
\begin{equation}\tilde{O}(\beta)=L(r,s)\#L(t,v)\label{222}
\end{equation}

where $\tilde{O}$ is irreducible, $\partial$-irreducible, and $d(\alpha,\beta)>1$. We now apply results on exceptional Dehn fillings at maximal distance. In particular, the following result of Boyer and Zhang \cite{Boyer} is relevant:

 \begin{thm}\label{boy}(\cite{Boyer}, Theorem 1.2.1) Let $M$ be a compact, connected, orientable, irreducible 3-manifold with $\partial M$ a torus. Assume that $M$ is neither a simple Seifert fibered manifold nor a cable on I(K) (the twisted $I$-bundle over the Klein bottle). Fix slopes $r_1$ and $r_2$ on $\partial M$ and suppose that $M(r_1)$ is a reducible manifold. If $M(r_2)$ has a cyclic fundamental group, then $d(r_1,r_2)\leq1$.
\end{thm}

In order to apply this result to Equations (\ref{111}) and (\ref{222}), we must check that $\tilde{O}$ satisfies the additional Seifert fibration and twisted $I$-bundle conditions. We say that a manifold is \textbf{simple} if it has no spheres, disks, tori or annuli that are essential.  We will proceed as follows: recall that a non-simple annular manifold with torus boundary is either a small Seifert fibered space or toroidal. We will show that $\tilde{O}$ cannot be a Seifert fibered space (SFS). Then, we will show that $\tilde{O}$ is atoroidal, which implies that it must be simple. This will show that $\tilde{O}$ is simple, which will allow us to use Theorem \ref{boy}.

\begin{prop}\label{prro} Let $M$ be a SFS such that $\partial M$ is a torus, and suppose that, for some slopes $\alpha$, $\beta$, $M(\alpha)=L(p,q)$  and $M(\beta)=L(r,s)\#L(t,u)$, with $r,t\neq1$. Then $d(\alpha,\beta)=1$.
\end{prop}

\begin{proof} We will consider two cases:
\begin{enumerate}
\item \textit{$r\neq2$ or $t\neq2$:} Suppose that $M$ admits a Seifert fibration. $M(\beta)$ is the result of performing Dehn filling on $M$ along a slope $\beta$ in $\partial M$. This filling always extends the fibration unless $\beta$ is isotopic to a regular fiber (see \cite{Brin}). Now we must distinguish two subcases: if both $r,t=0$, then $M(\beta)=S^2\times S^1\#S^2\times S^1$. This clearly does not admit a Seifert fibration. This implies that $\beta$ is the unique slope that destroys the fibration, and so $M(\alpha)$ extends such a fibration. This implies (see \cite{Ernst}) that the orbit surface of $M(\alpha)$ is either a torus or a nonorientable surface of genus two. This is impossible, since lens spaces can only fiber over $S^2$ or $\mathbb{R}P^2$ (here we are taking $L(0,1)$ as a lens space). If at least one of $r,t$ is not equal to zero, then $M$ either fibers over the disk with two exceptional fibers ($M$ cannot be a solid torus) or over $\mathbb{R}P^2$ with one exceptional fiber. In either case, attaching a solid torus as determined by $\alpha$ must add a regular fiber in order for us to obtain a lens space. Notice that once again $\beta$ destroys the fibration, so it must be isotopic to a fiber $H$ on $\partial M$. Choosing a crossing curve $Q$ for $\partial M$ and denoting a regular fiber by $H$, we must have that $\alpha\sim Q+\nu H$, and $\beta\sim H$, whence $d(\alpha,\beta)=1$. 

\item \textit{$r=t=2$}. Here $L(2,1)\#L(2,1)=\mathbb{R}P^3\#\mathbb{R}P^3$ can indeed be given a Seifert fibered structure. It is the orientable circle bundle over $\mathbb{R}P^2$, and thus it fibers over $\mathbb{R}P^2$ with no exceptional fibers (see \cite{Brin}). Therefore two situations may occur: either $M$ fibers over the disk with two exceptional fibers and we have obtained $L(2,1)\#L(2,1)$ by gluing a solid torus $V$ such that a meridian of $V$ is attached to a fiber $H$ on $\partial M$, or $M$ is the orientable circle bundle over the M\"obius band. If the former holds, then we proceed as in case (1) to conclude that $d(\alpha,\beta)=1$. If the latter holds, notice that there is an alternative fibering over the disk with two exceptional fibers, and hence this reduces to the previous situation, so we are finished.
\end{enumerate}
\end{proof}

Now, we show that $\tilde{O}$ is atoroidal. First, note that $\tilde{O}$ cannot contain non-separating tori, for these would remain non-separating after Dehn filling, contradicting the fact that all tori in lens spaces are separating. Next, we appeal to the following result, which can be found in \cite{Hat2}. This is a preliminary lemma in the proof of the JSJ decomposition Theorem. The result, we believe, is due to Kneser.

\begin{thm}\label{chop}Let $M$ be irreducible, compact and connected. Then, there exists a finite disjoint collection of tori $T$ such that each component of $M|T$ is atoroidal.
\end{thm}

Since all tori in $\tilde{O}$ are separating, the operation described above breaks $\tilde{O}$ into a disconnected collection of atoroidal manifolds. These are clearly, irreducible, $\partial$-irreducible. Since all boundary components are tori, each piece is either a small SFS or simple, as remarked before Propoposition \ref{prro}.

In order to aid our discussion, let us define the following:

\begin{defn} A \textbf{splitting graph of $\tilde{O}$ along $T$} is a graph $\mathcal{G}$ whose edges are the essential tori in $T$, and the vertices are connected pieces of the decomposition. Two vertices are connected by an edge if and only if some incompressible torus in $T$ separates their corresponding components, and $T$ lies in both components.  
\end{defn}

\begin{lemma} The splitting graph for $\tilde{O}$ is a tree.
\end{lemma}

\begin{proof} Recall that a finite connected graph is a tree if and only if removing any edge disconnects it. That $\mathcal{G}$ is connected follows from the fact that $\tilde{O}$ is connected. Since all of our tori are separating, deleting one of them will disconnect $\tilde{O}$. This implies that the splitting graph for this disconnected manifold is a disconnected graph, and the result follows.
\end{proof}

Now, we let the vertex corresponding to the connected component that contains $\partial \tilde{O}$ be the root of the splitting tree, and we can begin the analysis. Recall that in a tree, the concept of height of edges and vertices is well defined, once a root is chosen. We define the \textbf{level of a vertex} to be the number of edges that separate it from the root. We define the \textbf{level of an edge} to be the height of whichever of its vertices is closer to the root (where closer is to be understood in terms of levels). Also, note that there is only one vertex on the $0^{th}$ level: this is the root.

\vskip 0.1in We will number vertices as follows: $v_j^{k}$ will denote the $k^{th}$ vertex on the $j^{th}$ level, where the ordering at each level is chosen arbitrarily (we will soon see that this is irrelevant). An edge $e_j^k$ will be labelled similarly. We will abuse notation and identify vertices with pieces of the decomposition and edges with essential tori.

\vskip 0.1in Note that if $\mathcal{G}$ consists of a single vertex, there is nothing to prove: in this case, $\tilde{O}$ is atoroidal, and hence simple, by previous discussions. So we may assume that the set of edges is nonempty. We can now state our first result.

\begin{prop}\label{atoroidal}$v_0$ is a SFS over the annulus with at most one exceptional fiber.
\end{prop}

\begin{proof} As remarked above Proposition \ref{prro}, $v_0$ is either simple or a small SFS. Let us deal with the simple case first.
\par \textit{Case 1:} Suppose $v_0$ is simple. We may assume that $v_0$ has more than one boundary component: these are $\partial \tilde{O}$ and $\{e_0^{k}\}_{k=1}^{k=k_0}$, where $k_0$ is some positive integer. Now, perform the $\alpha$ and $\beta$ fillings along $\partial \tilde{O}$, and denote the manifolds obtained by $\alpha$ and $\beta$ filling on $v_0$ as $v_0(\alpha)$ and $v_0(\beta)$, respectively. We know that $L(p,q)$ and $L(r,s)\#L(t,v)$ are atoroidal, and hence all the $e_0^k$ must compress in $L(p,q)$ and $L(r,s)\#L(t,v)$.
\par We claim that at least one of the $e_0^k$ must compress in $v_0(\alpha)$, and the same for $v_0(\beta)$. Let us prove the claim: let $D$ be a compressing disk for  $e_0^1$ in $L(p,q)$. Since $e_0^1$ separates $\tilde{O}$, $int(D)$ is either contained in the component that contains $\partial \tilde{O}$ or in the other component. If $int(D)$ is in this other component, that means that $e_0^1$ compresses in $\tilde{O}$, contradicting the fact that it is essential. Thus we may assume $D$ is in the component which contains $\partial\tilde{O}$, but not necessarily in $v_0(\alpha)$. This will happen if and only if $D$ intersects some of the $e_0^k$. So suppose $D$ intersects some of the other $e_0^{k}$. This intersection is a finite collection of loops. Let $l$ be an innermost loop. We know that $l$ must then bound a disk $E\subset D$, and the boundary of this disk intersects, say, $e_0^m$. Notice that this disk doesn't intersect any other of the $e_0^k$. Now, we apply an innermost loop argument: if this disk is inessential in $e_0^m$, we may delete its interior from $e_0^m$, push $\partial E$ away from $e_0^m$, and glue back another disk disjoint from $e_0^m$, thus reducing the number of intersections. Hence, we can assume that $l$ is an essential loop in $e_0^m$ which bounds a disk. This is a compression disk for $e_0^m$. Repeating the argument at the start of the claim, we have that this disk is contained in the component to which $\partial \tilde{O}$ belongs. Since $E$ is taken not to intersect any of the other $e_0^k$, we have constructed a compressing disk for $e_0^m$ in $v_0(\alpha)$. Exactly the same argument goes for $v_0(\beta)$, and the claim is proved.
\par Now, we have a pair of fillings $\alpha$ and $\beta$ for a simple manifold $v_0$ such that $\Delta(\alpha,\beta)>1$ and such that $v_0(\alpha)$ and $v_0(\beta)$ are $\partial$-reducible. This contradicts results of Wu (Theorem 1 of \cite{Wu}), and so $v_0$ cannot be simple.
\par \textit{Case 2:} $v_0$ is a small SFS. We've already remarked that $v_0$ has at leat two boundary components. Going through the list of small SFS (see \cite{G2}), we find that the only possibilities are a SFS over the pair of pants with no exceptional fibers, or a SFS over the annulus with at most one exceptional fiber. We must rule out the pair of pants.
\par Suppose $v_0$ is a SFS over the pair of pants. We are given two different Dehn fillings on $v_0$. A filling with slope $\alpha$ ``destroys'' the fibration if and only if $\alpha$ is attached to a curve homologous to a fiber. Since we have two fillings with different slopes, at least one of them will ``respect" the Seifert fibration and WLOG we see that $v_0(\alpha)$ is a SFS over the annulus. We can now apply the same reasoning as in the simple case and deduce that at least one of the boundary components will compress in $v_0(\alpha)$. Now, it is known that an irreducible manifold with a compressible torus boundary component is a solid torus. Thus, $v_0(\alpha)$ is a solid torus, and this is clearly a contradiction. This completes the proof. \end{proof}

We must remark at this point that the case where $v_0$ is a SFS over the annulus with at most one exceptional fiber can occur: we can regard a SFS over the annulus with at most one exceptional fiber as the result of removing the neighborhood of a $(p,q)$ torus knot from a solid torus $V$, where this torus knot lies on the boundary of a torus concentric to $\partial V$. $q/p$ is the invariant of the exceptional fiber. Denote this space by $C_{p,q}$ (this is often known as a cable space). Performing Dehn filling on a boundary component yields a Seifert fibered space with at most two exceptional fibers for all slopes except one. Standard facts about Seifert fibered spaces can then be used to determine exactly when this fibration is that of a solid torus (that is, when no exceptional fibers are added). See \cite{Brin} for details.

\vskip 0.1in It might look that we cannot make any further progress with our analysis at this stage. However, we can get around the difficulties as follows: we claim that $v_0(\alpha)$ and $v_0(\beta)$ are both solid tori. Thus, we can think of $v_0(\alpha)$ and $v_0(\beta)$ as inducing two Dehn fillings on $v_1$, which we will denote by $v_1(\gamma)$, $v_1(\delta)$. If we can show that the distance between these fillings is not 1, we can apply Proposition \ref{atoroidal} inductively to show the tree is linear, with a terminal node $v_n$ which is a SFS over the annulus with at most one exceptional fiber and a pair of induced fillings on $v_n$ at distance greater than 1.

\begin{lemma} $v_0(\alpha)$ and $v_0(\beta)$ are solid tori.
\end{lemma}

\begin{proof} We know that at least one of the fillings will extend the Seifert fibration, so, WLOG, assume that $v_0(\alpha)$ is a SFS. In particular, it is irreducible. Exactly the same argument as for the pair of pants implies that $v_0(\alpha)$ is $\partial$-reducible, and hence a solid torus.
\par Now, if $v_0(\beta)$ also respects the fibration, there is nothing to prove. So assume that the fibration does not extend to $v_0(\beta)$. Choose a pair of curves $H$, $Q$ on the component $T$ of $\partial v_0$ where the filling takes place. Here, $H$ is a fiber, and $Q$ is a crossing curve. These form a basis for $H_1(T)$ (see \cite{Brin} for details on crossing curves). Thus, in the $\alpha$ filling, a meridian gets glued to a curve $\mu Q + \nu H$. In the $\beta$ filling, a meridian gets glued to $H$. However, we showed before that $v_0(\alpha)$ is a solid torus, and this cannot fiber with two exceptional fibers. This, together with the fact that $v_0$ has one exceptional fiber, implies that the $\alpha$ filling adds a regular fiber, and hence $\mu=1$. We thus have that $d(\alpha,\beta)=|1-0\cdot \nu|=1$. This contradicts $d(\alpha,\beta)>1$, and hence $v_0(\beta)$ must be a solid torus as well. Thus, the claim is proved.
\end{proof}

\vskip 0.1in
We must now compute the distance between these fillings $v_1(\gamma)$, $v_1(\delta)$. Essentially, this is figuring out which slopes on $\partial v_0(\alpha)$, $\partial v_0(\beta)$ bound disks. We use the following notation:  let $J$ be a simple closed curve in a solid torus $V$ which does not lie inside a ball, and let $Y=V-N(J)$, where $N(J)$ is a regular open neighborhood of $J$. We denote by $(J;r)$ the result of performing $r=m/n$ Dehn surgery on $Y$ on the boundary component created by removing $N(J)$.

The following result of Gordon is useful:

\begin{prop}(\cite{G4}, Lemma 3.3) Let $(J;r)$ and $V$ be as above. Furthermore, let $\mu$, $\lambda$ be a meridian/longitude basis for $\partial V$. Then, the kernel of the inclusion $i_*:H_1(\partial (J;r))\rightarrow H_1((J;r))$ is generated by:

$$
\left\{
\begin{array}{lr}
\frac{n\omega^2}{gcd(\omega,m)}\lambda+ \frac{m}{gcd(\omega,m)}\mu,\;\;\; \omega\neq 0\\
\mu,\;\;\;\;\;\;\;\;\;\;\;\;\;\;\;\;\;\;\;\;\;\;\;\;\;\;\;\;\;\;\;\;\;\; \omega=0
\end{array}
\right.
$$
\end{prop}

Here, $\omega$ is the winding number of $J$ in $V$. We can now use this to prove:

\begin{prop}$\Delta(\gamma,\delta)\geq2$.   
\end{prop}

\begin{proof}Let $\alpha=m/n$, $\beta=r/s$, with $\alpha$, $\beta$ as before. Then the previous proposition allows us to express $\gamma$ and $\delta$  in terms of coordinates on $\partial V$. We obtain that a meridian in $v_0(\alpha)$ is a $m/n\omega^2$ curve. Similarly, a meridian in $v_0(\beta)$ is a $r/s\omega^2$ curve. Hence, $\Delta(\gamma,\delta)=|\omega^2(ms-rn)|=\omega^2\Delta(\alpha,\beta)$. Since in the case of cabling spaces $\omega\neq 0$, the result follows.
\end{proof}

This allows us to complete our main argument:

\begin{thm}\label{mauro} $\tilde{O}$ is atoroidal.
\end{thm} 

\begin{proof} We've shown that $v_0$ is a SFS over the annulus, and that $v_0(\alpha)$ and $v_0(\beta)$ are solid tori. We can now pass the analysis to $v_1$, and so we have two Dehn fillings $v_1(\gamma)$ and $v_1(\delta)$. By the previous proposition, $\Delta(\gamma, \delta)\geq2$. We can now apply the argument of Proposition \ref{atoroidal} and conclude that $v_1$ is a SFS over the annulus with at most one exceptional fiber. Carrying out this procedure inductively we eventually reach the terminal node $v_n$ of our tree. We can easily see that $v_n$ cannot be simple, and it must therefore be a SFS over the disk with at most two exceptional fibers. But this is the situation of Proposition \ref{prro}: $v_n$ cannot be a SFS. We are forced to conclude that the collection of essential tori promised by Theorem \ref{chop} must be empty. That is, $\tilde{O}$ is atoroidal. 
\end{proof}

\begin{cor}\label{coro} $\tilde{O}$ is simple.
\end{cor}

\begin{proof} This is an immediate consequence of Proposition \ref{prro} and Theorem \ref{mauro}.
\end{proof}

And finally, we obtain:

\begin{cor}\label{coroo} There are no prime tangle solutions to Equations (\ref{tng1}), (\ref{tng2}), with $d(P,R)>1$. 
\end{cor}

\begin{proof} We've shown that $\tilde{O}$ is simple. Hence we can apply distance bounds on exceptional Dehn fillings. We have Dehn fillings $\tilde{O}(\alpha)$ yielding a lens space and $\tilde{O}(\beta)$ yielding a reducible manifold. This implies that $\Delta(\alpha,\beta)\leq1$, by Theorem \ref{boy}. The result follows.

\end{proof}

\subsubsection{\textit{Case $d(P,R)=1$}} Here, the Dehn filling techniques used above break down, and in fact the bounds obtained are best possible, in the sense that there are examples of reducible and cyclic fillings at distance one. The class of prime tangles (respectively irreducible, $\partial$-irreducible manifolds) is too large for one to hope for a complete solution to the problem. However, we can offer a complete solution for some families of tangles. We do this for the {Montesinos tangles}. This is done by classifying all solutions at distance one when $\tilde{O}$ is a Seifert fibred space. 

\begin{thm} Suppose that $\tilde{O}$ is a SFS, and that

\begin{equation} \tilde{O}(\alpha)=L(p,q)\label{a}
\end{equation}
\begin{equation} \tilde{O}(\beta)=L(r,s)\#L(t,u)\label{b}
\end{equation}

where $(r,s)$, $(t,u)$ are given pairs of integers. Then, $p=st+ru+rtm$ $(m\in\mathbb{Z})$, and $q=sk+ry$, where $k,y$ are any integers such that $k(u+tm)-yt=1$.
\end{thm}

\begin{proof} We've already remarked that the case $r=t=0$ cannot happen. Suppose that both $r$ and $t$ are not equal to zero. Then, by the remarks in Proposition \ref{prro}, $\widetilde{O}$ fibers over the disk with exactly two exceptional fibers, of type $r/s$ and $t/u$ (see \cite{Brin}). These spaces are double branch covers of the {Montesinos tangles}. We can exploit this fact to simplify our computation. The Seifert fibred space over the disk with two exceptional fibers of type $r/s$, $t/u$ is a double branched cover of the tangle $O=T(\frac{r}{s})+T(\frac{u}{t})$. It is a straightforward matter to verify that $N(O+T(\infty))$ is a nontrivial connected sum. This is the only tangle addition that has this property, and it corresponds to adding a fibered solid torus such that the meridian gets glued to a fiber $H$ in $\partial \widetilde{O}$. On the other hand, the fillings that give lens spaces correspond to the integral tangle additions downstairs. Using the facts that $T(\frac{u}{t})+T(m)=T(\frac{u+tm}{t})$ and $N(\frac{a}{b}+\frac{c}{d})=N(\frac{ad+bc}{ad'+bc'})$, for any $c'$ and $d'$ such that $c'd-cd'=1$ (see \cite{Sum}) we have that $N(\frac{s}{r}+\frac{u+tm}{t})=N(\frac{st+ru+rtm}{sk+ry})$, where $k,y$ are such that $k(u+tm)-yt=1$. This gives all the possibilities as one varies $m$.
\end{proof}

This result can easily be translated into the language of Montesinos tangles.

\begin{cor} With $p=r=t$, Equations (\ref{a}), (\ref{b}) are impossible if $\widetilde{O}$ is a Seifert fibered space.
\end{cor}

\begin{proof} If $p=r=t$, then we must have $2p+p^2m=p$, whence $pm=-1$. Since we are assuming that $|p|>1$, we have a contradiction.
\end{proof}

\vskip 0.1in

\section{Solution for processive recombination}

Here, we apply the results of Ernst and Sumners \cite{Sum} to equations \ref{1}-\ref{4}, corresponding to processive recombination of Hin on mutated sites. Although these results have not been published, they are known to the community, and we display them here to illustrate the tangle model. 

\vskip 0.1in\textbf{Note:} We will not solve here the tangle equations for Hin acting processively on standard sites. Biologically, this is completely analogous to the action of the protein Gin, characterised by Vazquez and Sumners \cite{Vaz}. It yielded the solution $O=T(-1/2), R=(1)$.

\vskip 0.1in As discussed in Section 3, we show that the tangles involved must be rational. A straightforward observation based on the unique factorisation of knots implies that $O$, $P$ and $R$ cannot be locally knotted. Thus we are left to rule out prime tangles.

\vskip 0.1in  The proof that $R$ must be rational is identical to Lemma 3.1 in \cite{Sum}. In the following, we show that $O$ must be rational, as well. 

\begin{prop}\label{same}$O$ is rational.
\end{prop}

\begin{proof}We know $O$ is not locally knotted, so suppose $O$ is prime. Then, it can be shown that $P$ is rational (\cite{Sum}, Lemma 3.1). This implies that $\tilde{O}$ is the complement of a non-trivial knot $S^3$ (see \cite{Bl}). Suppose that it is the complement of a non-torus knot, say $\tilde{O}=M_{K}$. Equations (\ref{2}) and (\ref{3}) imply, in the double branched cover setting, $\tilde{O}\cup_g V_1\cong L(3,1)$ and  $\tilde{O}\cup_h V_2\cong L(7,3)$, \textit{i.e.} we have two Dehn surgeries on $K$ yielding lens spaces. By the Cyclic Surgery Theorem (see \cite{CGLS}), they must be integral. Hence, $\tilde{O}(m)\cong L(3,1)$ and $\tilde{O}(n)\cong L(7,3)$, for some integers, $m,n$. An easy homology argument implies that $m=3, n=7$. But this contradicts the fact that the integers must be successive. Hence we may assume that $K$ is a  $(r,s)$ non trivial torus knot. Moser \cite{Mos} shows that, for $p/q$ surgery on a $(r,s)$ torus knot, we obtain a lens space if and only if $|rsq+p|=1$, and the resulting lens space is $L(p,qs^2)$. Now, $r,s$ are positive integers, and we have $|p|=rs|q|\pm1$. Since we have $L(3,1)$ as a result of surgery, we wish to solve $|3|=rs|q|\pm1$. But for a torus knot to be non-trivial, we must have $r\geq3$, $s\geq2$, and hence the above equation cannot have integer solutions. We conclude that $K$ must be the unknot. The result follows.
\end{proof}

We now calculate the tangle solutions. In the following discussion, we warn the reader that we will not take chirality into account, and so we will take liberties such as to consider $b(3,1)\sim b(3,-1)$.

\begin{thm}Given Equations (\ref{1})-(\ref{4}), then either $O=T(1/2)$ and $R=T(-2)$, or $O=T(-1/2)$ and $R=T(2)$.
\end{thm} 

\begin{proof} This boils down to a computation (\cite{Sum}, Lemma 2.1). We first check that $R$ must be integral. Suppose not. Then, the double branched cover of $N(O+R+R+R)$ is a Seifert fibered space over $S^2$ with at least 3 exceptional fibers. Now, all prime 7-crossing knots are 4-plats (see the table at the end of \cite{Burde}), and we know that lens spaces cannot fiber with more than two exceptional fibers. Hence we obtain a contradiction, so $R$ is integral.
\par Let $R=T(r)$. Clearly, $R+R=T(2r)$. Let $O=T(u/v)$. Applying Lemma 2.1 from \cite{Sum}, we obtain the following equations $|u+k|=3$ and $|u+2k|=7$, where $k=rv$. Solving gives the following solution pairs: $(u,k)=\{(1,-4),(-1,4),(2,-5),(-2,5)\}$. We must now check case by case.
\par Suppose $r=1, v=-4$ (the reverse sign case is essentially the same, but with mirror images). We have, upon checking, $N(O+R)=b(3,1)$, $N(O+R+R)=b(7,3)$, which is consistent. However, $N(O+R+R+R)=b(11,3)\sim b(11,4)$. From tables \cite{Burde}, this is Stevedore's knot, a 6-crossing knot. This is inconsistent, so we discard this solution,
\par Suppose $r=4, v=-1$. We obtain $N(O+R+R)=b(7,8)\sim b(7,1)$, which is inconsistent, so we discard this solution.
\par Suppose $r=2, v=-2$. We obtain $N(O+R)=b(3,1)$, $N(O+R+R)=b(7,3)$, $N(O+R+R+R)=b(11,6)\sim b(11,-5)\sim b(11,5)$. This is consistent, so we have found a solution ($r=-2, v=2$ is also a solution, and this gives mirror images). The other cases are checked similarly, and all the solutions are discarded.
\end{proof}

\begin{cor}Under the hypotheses of the tangle model, $N(O+R+R+R)=b(11,5)$.
\end{cor}

\begin{proof} This is a straightforward computation.
\end{proof}

We've isolated the unique tangle solution, modulo mirror images. We can now decide on the chirality by looking again at Figure \ref{lhtrefoil}. It is clear that the trefoil obtained is a left-handed one. As only the pair $O=T(-1/2)$, $R=T(2)$ yields  $N(O+R)=b(3,1)$ as the left handed trefoil,  we take it as the solution. Notice that this determines the handedness of all subsequent products of processive recombination. Projections of the solution tangles are displayed in Figure 9.

\begin{figure}[h!]
\centering
\includegraphics[scale=0.5]{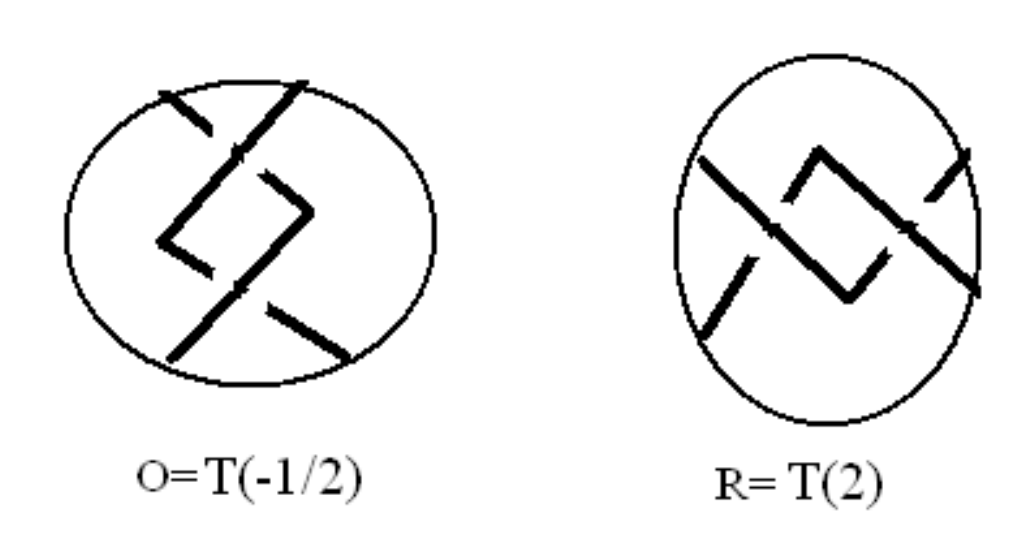}
\caption{Solutions for processive recombination.}

\end{figure}

Unfortunately, little can be said about the tangle $P$. In fact, it admits infinitely many solutions (see \cite{Sum}). This is essentially due to the fact that it only appears once in the equations, and this does not give enough information to determine its type. Following Vazquez (see \cite{Vaz} and references therein), we adopt the common approach and set $P=T(0)$.

\begin{figure}[h!]
\centering
\includegraphics[scale=0.4]{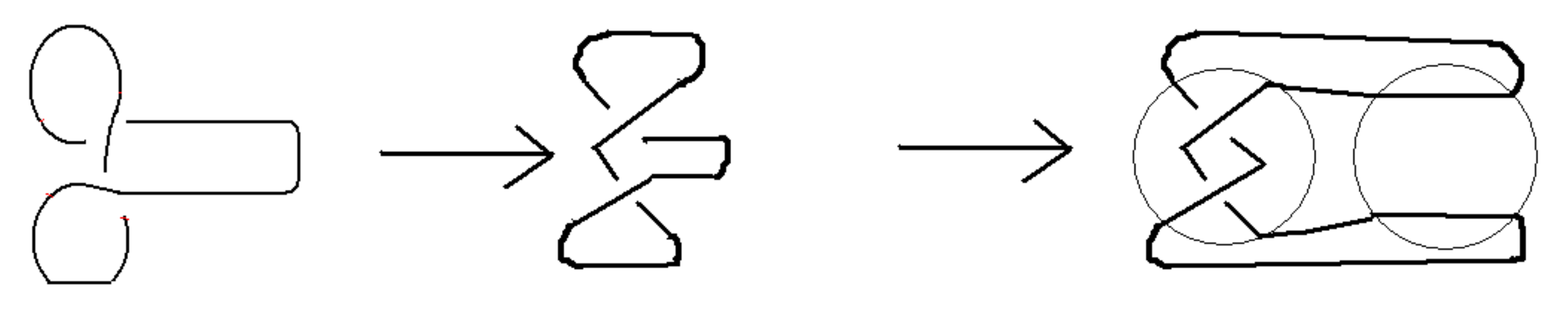}
\caption{Diagram of an isotopy of the pre-recombination DNA into $N(T(-1/2)+T(0))$.} 
\label{pickw}
\end{figure}

The results obtained confirm the model of Heichman \textit{et al} \cite{Heichman} that  recombination of the mutated sites occurs by a 360 degree  clockwise twist at the recombinations sites, represented by the tangle $T(2)$. The $T(-1/2)$ solution is consistent with the fact that the enhancer sequence does indeed acquire the predicted relative position with respect to the crossover sites. This is in agreement with the treatment of Vazquez \cite{Vaz} The solutions are schematised in Figure \ref{pickw}.

\section{Solution for distributive recombination}

We now tackle the equations for distributive recombination. Recall from Section 3 that these were

$$N(Q+P)=b(3,1)
$$

$$N(Q+R)=b(3,1)\#b(3,1)
$$

 This will characterise, topologically, all the configurations of the outside bound DNA during distributive recombination steps. We take $P=T(0)$ and $R=T(2)$, as determined in the previous section. This is essentially the assumption that the mechanism of action of the protein is independent of whether recombination is processive or distributive.  This reflects the experimental evidence \cite{reid}.

To prove Theorem 6.2, we need the following result of Ernst and Sumners:

\begin{lemma}\label{lil}(Lemma 3.6 \cite{Sum}) Let $A=T(0)$ and $B=T(\alpha/\beta)$, and suppose $X$ is a tangle such that $N(X+A)=b(1,1)$ and $N(X+B)=b(p,q)$. Then, $L(p,q)$ can be obtained by $(\beta+s\alpha)/\alpha$ Dehn surgery on $M_{K}=\tilde{X}$, where $s\in \mathbb{Z}$.
\end{lemma}

\begin{thm}\label{ti}Given that $P=T(0)$, $R=T(2)$, the following are the complete tangle solutions for $Q$ to the equations:

$$N(Q+P)=b(3,1)
$$
$$N(Q+R)=b(3,1)\#b(3,1)
$$
 
\begin{enumerate}
\item \textbf{Rational tangles:} There are no rational tangle solutions.
\item \textbf{Prime tangles:} There are no prime tangle solutions.
\item \textbf{Locally knotted solutions:} The only 4 possible solutions for $Q$ are displayed in Figure \ref{knotty}.
\end{enumerate}

\begin{figure}[h!]
\centering
\includegraphics[scale=0.4]{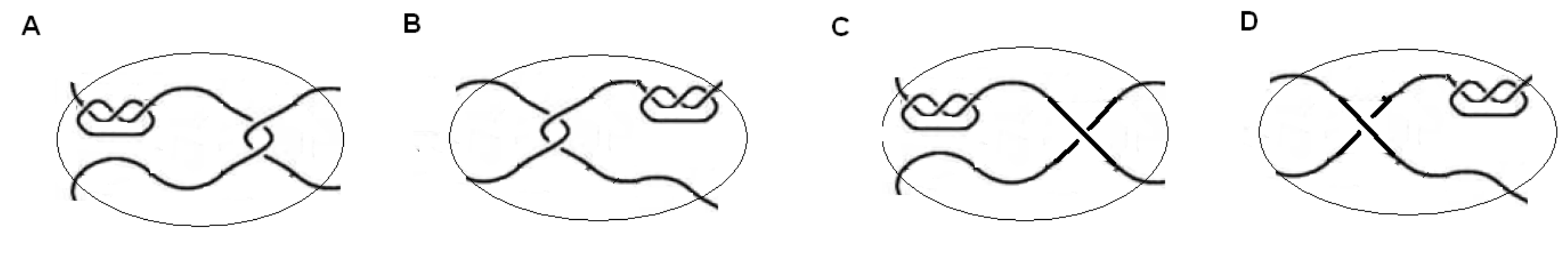}
\caption{The four topological possibilities for $Q$. Note that only A and B are consistent with experimental evidence.}
\label{knotty}
\end{figure}
\label{final}

\end{thm}

\begin{proof} We treat each case separately:

\par $(1)$ This is a straightforward observation, already remarked in Section 4.
\par $(2)$ This follows from Corollary \ref{coroo}.
\par $(3)$  We use the procedure described at the beginning of Section 4.1 to obtain a one-to-one correspondence between locally knotted solutions and solutions to the equations obtained by excising a common knotted (ball, arc) pair.

$$N(\hat{Q}+P)=b(1,1)
$$
$$N(\hat{Q}+R)=b(3,1)
$$

Here, $\hat{Q}$ is either rational or prime. Suppose $\hat{Q}$ is prime. Then, its double branched cover is a knot complement $M_{K}$. Suppose $K$ is not a torus knot (\textit{i.e.} $\tilde{Q}$ is not a SFS). By Lemma \ref{lil}, we can obtain $L(3,1)$ by performing $(1+2s)/2$ surgery on $K$, contradicting the Cyclic Surgery Theorem (which asserts that the surgery must be integral). Hence we may assume that $K$ is a torus knot. Applying exactly the same argument as in Proposition \ref{same}, we conclude that $K$ is the unknot, \textit{i.e.} $\hat{Q}$ is rational.
\par It is now a routine matter to enumerate the possibilities for $\hat{Q}$. Let $\hat{Q}=T(u/v)$. Applying Lemma 2.1 of \cite{Sum} to both equations with $P=T(0)$ and $R=T(2)$ we get $|u|=1$ and $|u+2v|=3$. Solving gives the rational numbers $u/v=\{1/1, -1/2\}$. Reintroducing the excised trefoils, we obtain the desired locally knotted solutions.

\end{proof}

The results suggest that, with the exception of solutions C and D in Figure \ref{knotty}, the action of Hin during distributive recombination is essentially identical to that of processive recombination. In particular, this seems to indicate the formation of composite knots occurs as Hin traps a DNA molecule which has already undergone one round of recombination, with the existing DNA knot being localised outside of the invertasome. The invertasome maintains a configuration that is identical to the processive recombination steps. This is consistent with the fact that Hin is a highly restrictive protein, and the invertasome tends to be productive only under very special conditions \cite{reid}. These results are what one would naturally expect of distributive recombination yielding composite knots.
\par The undesired solutions illustrated in Figures \ref{knotty}C and D do not seem to conform to any previously reported biological mechanism \cite{reid}. In particular, they induce an antiparallel alignment of the inverted crossover sites. However, as discussed in Section 2, the DNA segments involved in the invertasome are only several dozen bp long. This short length makes it highly unlikely a compensating supercoil could be introduced within the invertasome to give the required parallel alignment of the inverted crossover sites.

\section{Conclusions and directions for future research}

We have successfully extended the tangle model to treat equations involving composite knots in most cases. To do this, we first built on results on exceptional Dehn fillings at maximal distance. We then showed that, for any prime tangle, there are no rational tangle attachments of distance greater than one that first yield a 4-plat and then a connected sum of 4-plats. Next, we applied the tangle model to describe the rational tangle solution to processive Hin recombination on mutated sites. Finally, this was combined with our surgery results to classify all tangle solutions arising from distributive recombination mediated by Hin (Theorem \ref{ti}).
\par Our work confirms the biologically predicted configuration for the invertasome (modulo one undesired solution, that can be excluded on biological grounds). Further, the outside DNA attains a configuration that is the expected one, with localised knotted arcs trapped outside the invertasome, as shown in Figure \ref{knotty}. 
\par We emphasize that although we have illustrated this extended tangle model with the Hin recombinase, that the lens space connect sum surgery restrictions, as developed in Section 4, are completely general. Thus, these can be applied to a variety of site-specific recombination systems where composite knots arise as a result of distributive recombination. In particular, the Gin recombinase has been known to yield composite knots \cite{kan}.   

\subsection{Directions for future research} Experimentally, it would be useful to obtain electron micrograph images of reaction products of the form $b(3,1)\#b(7,3)$, which we predict would occur as the result of an additional round of recombination in our model. The existence of these products would eliminate solutions C and D in Theorem \ref{ti}, as discussed in the previous section. 
\par As remarked in Section 3, experiments indicate the Hin invertasome behaves as a 3-string tangle. One other possible avenue of research would therefore be to model Hin as a 3-string tangle. Darcy, Luecke and Vazquez \cite{Darcy} and Cabrera Ibarra \cite{Ib}, \cite{Ib2}, \cite{Ib3} have published some results in this direction. In particular, Cabrera Ibarra gives results on the classification of 3-string tangle solutions to processive recombination by Gin, a protein very similar to Hin. We plan to extend these results to the setting of distributive recombination.
\par Finally, the most natural direction in which to extend our results would be to obtain more information on the $d(P,R)=1$ case. This could be then be applied to a variety of site-specific recombination systems, including distributive recombination of Hin on standard sites \cite{Heichman}, and processive followed by distributive recombination of Gin on directly repeated sites \cite{Sum3}, \cite{kan}. Additionally, classifying such surgeries at distance one would be an interesting mathematical problem on its own. We are currently using knot Floer homology \cite{KMOS} to address this question.

\section*{Acknowledgements} We wish to thank, Andrew Lobb, Koya Shimokawa and DeWitt Sumners for insightful discussions. We are particularly indebted to Ken Baker for invaluable input, and for commenting on an early draft.
DB is supported in part by EPSRC Grants EP/H0313671, EP/G0395851 and  EP/J1075308, as well as an LMS Scheme 2 Grant. MM is supported by an EPSRC DTA.

 \bibliographystyle{amsalpha}

\begin{thebibliography}{A}

\bibitem{Bl}S. Bleiler
\textit{Prime tangles and composite knots}, Knot theory and manifolds. Proceedings 1983, Springer-Verlag, 1985

\bibitem{Boyer}S. Boyer, X. Zhang, \textit{On Culler-Shalen seminorms and Dehn filling}, Annals of Math., 148, 737-801, 1998

\bibitem{Boyer2}S. Boyer, \textit{Dehn Surgery on Knots}, Chaos, Solitons and Fractals, vol. 9, 657-670, 1998 

\bibitem{Brin}M. Brin,
\textit{Seifert Fibered Spaces (lecture notes)},
ftp://ftp.math.binghamton.edu/pub/matt/seifert.pdf ,
1993

\bibitem{Buck2} D. Buck and C. Verjosvky Marcotte, \textit{Tangle solutions for a family of DNA rearranging proteins}, Math. Proc. Cambridge Philos. Soc., 139, n0.1, 59-80, 2005

\bibitem{Buck} D. Buck and C. Verjosvky Marcotte, \textit{Classification of tangle solutions for integrases, a protein family that changes DNA topology}, Journal of Knot Theory and its Ramifications, 16, 2007

\bibitem{Burde}G. Burde, H. Zieschang, \textit{Knots}, 
de Gruyter Studies in Mathematics, 2003

\bibitem{Ib}H. Cabrera Ibarra, \textit{On the classification of rational 3-tangles}, Journal of knot theory and its ramifications, Vol. 12, 7, 921-946, 2003

\bibitem{Ib2} H. Cabrera Ibarra, \textit{Braid solutions to the action of the Gin enzyme}, Journal of Knot Theory and its ramifications (in press)

\bibitem{Ib3} H. Cabrera Ibarra, \textit{An algorithm based on 3-braids to solve tangle equations arising in the action of Gin DNA invertase}, Applied Mathematics and Computation, Vol. 216, Issue 1, 95-106, 2010 

\bibitem{Conway} J. Conway, \textit{An enumeration of knots and links, and some of their algebraic properties}, in Computational Problems in Abstract Algebra, 329-358, Pergamon Press, 1969

\bibitem{Cozz}N. Cozzarelli, S. Wasserman, \textit{Biochemical Topology: applications to DNA recombination and replication}, Science, vol. 232, Issue 4753, 951-960, 1986

\bibitem{Cris}N. Crisona R. Weinberg, B. Peter. D.W. Sumners, \textit{The topological mechanism of phage lambda integrase}, J. Mol. Biology, 289, no. 4, 747-775, 1999

\bibitem{CGLS} M. Culler, C. McA. Gordon, J. Luecke, P. Shalen \textit{Dehn Surgery on Knots}, Bulletin of the AMS, 1985

\bibitem{Darcy2}I. Darcy, \textit{Biological distances on DNA knots and links: applications to XER recombination}, Journal of knot theory and its ramifications, no.2, 269-294, 2001

\bibitem{Darcy}I. Darcy, J. Luecke, M. Vazquez, \textit{Tangle analysis of difference topology experiments: applications to a Mu protein-DNA complex}, Algebraic and Geometric Topology 9, 2247-2309, 2009 

\bibitem{Ernst} C. Ernst, \textit{PhD dissertation}, Florida State University, 1989

\bibitem{Sum}C. Ernst, D.W. Sumners,
\textit{A calculus for rational tangles: applications to DNA recombination},
Math. Proc. Camb. Phil. Soc., 108, 409,
1990


\bibitem{Ernst2} C. Ernst, \textit{Solving tangle equations}, Journal of Knot Theory and its ramifications, no. 2, 145-159, 1996

\bibitem{Ernst3} C. Ernst, \textit{Solving tangle equations II}, Journal of Knot Theory and its ramifications, no. 1, 1-11, 1997


\bibitem{G4} C. McA. Gordon and J. Luecke, \textit{Reducible manifolds and Dehn Surgery}, Topology, 35, 385-409, 1996

\bibitem{G2}C. McA. Gordon,
\textit{Dehn Surgery},
unpublished course notes,
2006

\bibitem{G3}C. McA. Gordon, 
\textit{Small surfaces and Dehn Filling}, 
Geometry and Topology Monographs, vol. 2, 177-199
1998

\bibitem{G4}C. McA. Gordon, \textit{Dehn Surgery and satellite knots}, Trans. of the AMS, vol.275, 2, 1983

\bibitem{Grindley}N. Grindley, K. Whiteson, P. Rice 
\textit{Mechanisms of Site-Specific Recombination}, 
Annu. Rev. Biochem., 75, 567-605, 
2006

\bibitem{Hat2} A. Hatcher,
\textit{Notes on Basic 3-Manifold Topology},
http://www.math.cornell.edu/$\sim$ hatcher/3M/3Mdownloads.html,
2003

\bibitem{Heichman}K. Heichman \textit{et al}, \textit{Configuration of DNA strands and mechanism of strand exchange in the Hin invertasome as revealed by analysis of recombinant knots}, Genes and Development, 5, 1622-1634, 1991

\bibitem{Hodge}C. Hodgson and J. Rubinstein, \textit{Involutions and isotopies of lens spaces}, Knot Theory and manifolds, proc., Vancouver, 1983

\bibitem{reid}R. Johnson, \textit{Bacterial site-specific DNA inversion systems}, in \textit{Mobile DNA II}, edited by N. Craig \textit{et al}, 2002 AMS press

\bibitem{kan}R. Kanaar \textit{et al}, \textit{Processive recombination by the phage Mu Gin system: implications for the mechanisms of DNA strand exchange}, Proc. natn. Acad. Sci. USA, 85, 752-756, 1990

\bibitem{Katsu}K. Kutsukake \textit{et al}, \textit{Two DNA Invertases Contribute to Flagellar Phase Variation in Salmonella enterica Serovar Typhimurium Strain LT2}, Journal of Bacteriology, vol. 188, no. 3, 950-957, 2006

\bibitem{Kim}P. K. Kim, J. Tollefson, \textit{Splitting the PL involutions of non-prime 3-manifolds}, Michigan Math. J., 27, 259-274, 1980

\bibitem{KMOS}P. Kronheimer, T. Mrowka, P. Ozsv\'ath, Z. Szab\'o, \textit{Monopoles and lens space surgeries}, Ann. Math., 165, 457-546, 2004.

\bibitem{Lick}R. Lickorish, 
\textit{Prime knots and tangles}, Transactions of the AMS, vol. 267, 1981

\bibitem{Monte} J. Montesinos, \textit{Variedades de Seifert que son recubridores ciclicos ramificados de dos hojas}, Boletin de la Sociedad Matematica Mexicana, vol. 18, 1-32, 1973.

\bibitem{Mos}L. Moser,
\textit{Elementary surgery on a torus knot},
Pacific Journal of Mathematics, vol. 38, no. 3,
1971

\bibitem{Rol} D. Rolfsen,
\textit{Knots and Links},
Publish or Perish,
1975

\bibitem{Sum3}D.W. Sumners, C. Ernst, S. Spengler, N. Cozzarelli \textit{Analysis of the mechanism of DNA recombination using tangles}, Quart. Rev. Biophysics, 254-313, 1995

\bibitem{Sum2}D.W. Sumners, \textit{Lifting the curtain: using topology to probe the hidden action of protein}, Notices of the AMS, vol. 42, no.5, 1995.


\bibitem{Vaz}M. Vazquez and  D.W. Sumners, \textit{Tangle analysis of Gin site-specific recombination}, Math. Proc. Camb. Phil. Soc., 136, 2004

\bibitem{Vaz2}M, Vazquez, S.D. Colloms and D.W. Sumners, \textit{Tangle analysis of Xer recombination yields only three solutions, all consistent with a single three dimensional topological pathway}, J. Mol. Biology, 346, no. 2, 493-504, 2005

\bibitem{Wu}Y. Q. Wu, \textit{Incompressibility of surfaces in surgered 3-manifolds}, Topology 31, 271-289, 1992

\end{thebibliography}

\end{document}